\documentstyle[12pt]{article}
\begin{document}
\input{amssym}
\title{Exact solution of generalized \\ inviscid Burgers' equation}
\author{Mehdi Nadjafikhah\thanks{Faculty of Mathematics, Iran University of
Science and Technology, Narmak, Tehran, I.R.IRAN. e-mail:
m\_nadjafikhah@iust.ac.ir} }
\date{}
\maketitle
\begin{abstract}
Let $f,g:{\Bbb  R}\to{\Bbb R}$ be integrable functions, $f$
nowhere zero, and $\varphi (u)=\int du/f(u)$ be invertible. An
exact solution to the generalized nonhomogeneous inviscid Burgers'
equation $u_t+g(u).u_x=f(u)$ is given, by quadratures.
\end{abstract}
{\bf 2000 Mathematics subject classification:} 35Q53, 37K10.\\
{\bf Keywords.} KdV-like equations, Burgers' equation.
\input{amssym}
\def\be{\begin{eqnarray}}
\def\ee{\end{eqnarray}}
\def\di{\displaystyle}
\def\rank{{\bf rank}}

\section{Introduction}
Burgers' equation is a fundamental partial differential equation
from fluid mechanics (see, \cite{[5]}, \cite{[1]} and \cite{[2]}).
It occurs in various areas of applied mathematics, such as
modeling of gas dynamics and traffic flow. It is named for
Johannes Martinus Burgers (1895-1981). For a given velocity $u$
and viscosity coefficient $\nu$, the general form of Burgers'
equation is: $u_t+u.u_x=\nu.u_{xx}.$ When $\nu=0$, Burgers'
equation becomes the inviscid Burgers' equation: $u_t+u.u_x=0,$
which is a prototype for equations for which the solution can
develop discontinuities (shock waves).

The generalized nonhomogeneous inviscid Burgers' equation
\cite{[4]} is $${\cal E}_{u,g,f}: u_t+g(u).u_x=f(u).$$ We solve
this equation by quadratures; at two stages depending only on the
function $f$.
\section{Homogeneous Burgers' equation}
\paragraph{2.1 Theorem.}
{\it Let $g:{\Bbb R}\to{\Bbb R}$ be an arbitrary function. Then,
$$u=h(x-t.g(u))$$ gives the solution to ${\cal
E}_{u,g,0}:u_t+g(u).u_x=0$, where $h:{\Bbb R}\to{\Bbb R}$ is an
arbitrary function.}

\bigskip \noindent {\it Proof:} Let $\tilde{x}=x-t.g(u)$,
$\tilde{t}=t$ and $\tilde{u}(\tilde{x},\tilde{t})=u$ are new
independent and dependent variables. Then,
$x=\tilde{x}+\tilde{t}.g(\tilde{u})$, $t=\tilde{t}$, and
$u=\tilde{u}$. Now, $u_{\tilde{x}}=\tilde{u}_{\tilde{x}}$ and
$u_{\tilde{x}} = \big(u(x,t)\big)_{\tilde{x}} =
u_x.x_{\tilde{x}}+u_t.t_{\tilde{x}} =
u_x.(1+\tilde{t}.g'(\tilde{u}).\tilde{u}_{\tilde{x}})$. Therefore
\be
u_x.(1+\tilde{t}.g'(\tilde{u}).\tilde{u}_{\tilde{x}})=\tilde{u}_{\tilde{x}}.
\label{eq:3} \ee
Similarly, $u_{\tilde{t}}=\tilde{u}_{\tilde{t}}$ and
$u_{\tilde{t}} = \big(u(x,t)\big)_{\tilde{t}} =
u_x.x_{\tilde{t}}+u_t.t_{\tilde{t}} =
u_x.(g(\tilde{u})+\tilde{t}.g'(\tilde{u}).\tilde{u}_{\tilde{t}})+u_t$.
Therefore
\be
u_x.(g(\tilde{u})+\tilde{t}.g'(\tilde{u}).\tilde{u}_{\tilde{t}})+u_t=\tilde{u}_{\tilde{t}}.
\label{eq:4}\ee
By solving (\ref{eq:3}) and (\ref{eq:4}), we have $u_x =
\tilde{u}_{\tilde{x}}/\Delta$, $u_t
=(\tilde{u}_{\tilde{t}}-g(\tilde{u}).\tilde{u}_{\tilde{x}})/\Delta$,
where $\Delta=1+\tilde{t}.g'(\tilde{u}).\tilde{u}_{\tilde{x}}$.
Thus, the equation $u_t+g(u).u_x=0$ reduced to
$\tilde{u}_{\tilde{t}}=0$. Therefore, $\tilde{u}=h(\tilde{x})$,
and Theorem was proved. \hfill\ $\Box$
\paragraph{2.2 Example.}
As an illustration of Theorem 2.1, consider the equation
\be u_t+u.u_x=0. \label{eq:p1}\ee
Then $g=u$, and by the Theorem 2.1, the solution of this equation
is $u=h(x-t.u)$, where $h:{\Bbb R}\to{\Bbb R}$ is an arbitrary
function. We look at a couple of special cases:
\begin{itemize}
\item If $h(s)=(a.s+b)/c$, where $a$, $b$ and $c\neq0$ are
arbitrary constants; Then $$\di u(x,t)=\frac{a.x+b}{a.t+c}$$ is a
solution to the (\ref{eq:p1}).
\item If $h(s)=\sqrt{a.s+b}+c$, where $a$, $b$ and $c$ are
arbitrary constants; Then
$$\di
u(x,t)=c-\frac{1}{2}.\Big(a.t\pm\sqrt{a^2.t^2+4.a.x-4.a.c.t+4.b}\Big)$$
is a solution to the (\ref{eq:p1}).
\item If $h(s)=e^{a-s}$, where $a$ is
an arbitrary constant; Then
$$\di u(x,t)=-\frac{1}{t}.LW(-te^{a-x})$$ is a solution to the
(\ref{eq:p1}), where $LW$ is the Lambert $W-$function; i.e. a
function defined by function-equation $f(x).e^{f(x)}=x$.
\end{itemize}
\section{Nonhomogeneous Burgers' equation}
\paragraph{2.1 Lemma.}
{\it Let $f:{\Bbb R}\to{\Bbb R}$ be a nowhere zero continuous
function, and $v=\varphi (u)=\int du/f(u)$. Then, every ${\cal
E}_{u,g,f}$ can be transformed into ${\cal E}_{v,g\circ \varphi
^{-1},1}$.}

\bigskip \noindent {\it Proof:}
By continuity of $f$, $\varphi$ is a smooth invertible function,
and
\begin{eqnarray*} v_t+(g\circ \varphi ^{-1})(v).v_x =\varphi _t(u)+g(u).\varphi _x(u)
= \frac{1}{f(u)}.\Big(u_t+g(u).u_x\Big)= 1.
 \end{eqnarray*}
Therefore, $v$ satisfy in a ${\cal E}_{v,g\circ \varphi ^{-1},1}$
equation. \hfill\ $\Box$
\paragraph{3.2 Lemma.}
{\it Let $g$ be an integrable function. Then, the solution of
${\cal E}_{u,g,1}$ is $x=\int g(u)\,du+h(t-u)$, where $h:{\Bbb
R}\to{\Bbb R}$ is a arbitrary function.}

\bigskip \noindent {\it Proof:}
Let $\tilde{x}=u$, $\tilde{t}=t-u$ and
$\tilde{u}(\tilde{x},\tilde{t})=x$ are new independent and
dependent variables. Then, $x=\tilde{u}$, $t=\tilde{t}+\tilde{x}$,
and $u=\tilde{x}$. Now, $u_{\tilde{x}}=(\tilde{x})_{\tilde{x}}=1$
and $u_{\tilde{x}} = \big(u(x,t)\big)_{\tilde{x}} =
u_x.x_{\tilde{x}}+u_t.t_{\tilde{x}} =
u_x.\tilde{u}_{\tilde{x}}+u_t$. Therefore
\be u_x.\tilde{u}_{\tilde{x}}+u_t=1. \label{eq:1} \ee
Similarly, $u_{\tilde{t}}=(\tilde{x})_{\tilde{t}}=0$ and
$u_{\tilde{t}} = \big(u(x,t)\big)_{\tilde{t}} =
u_x.x_{\tilde{t}}+u_t.t_{\tilde{t}} =
u_x.\tilde{u}_{\tilde{t}}+u_t$. Therefore
\be u_x.\tilde{u}_{\tilde{x}}+u_t=0. \label{eq:2}\ee
By solving (\ref{eq:1}) and (\ref{eq:2}), we have $u_x =
1/(\tilde{u}_{\tilde{x}}-\tilde{u}_{\tilde{t}})$, $u_t =
-\tilde{u}_{\tilde{t}}/(\tilde{u}_{\tilde{x}}-\tilde{u}_{\tilde{t}})$,
and equation $u_t+g(u).u_x=1$ reduced to
$\tilde{u}_{\tilde{x}}=g(\tilde{x})$. Therefore, $\tilde{u}=\int
g(\tilde{x})\,d\tilde{x} + h(\tilde{t})$. This proves Lemma.
\hfill\ $\Box$

\bigskip The proof of following Theorem relies on Lemmas 3.1 and
3.2.
\paragraph{3.3 Theorem.}
{\it Let $f,g:{\Bbb R}\to{\Bbb R}$ be integrable functions, $f$
nowhere zero, $ \varphi (u)=\int du/f(u)$ invertible, and
$\ell(u)=\int (g\circ \varphi ^{-1})(u)\,du$. Then, $$x=(\ell\circ
\varphi )(u)+h\big(t-\varphi (u)\big),$$ gives the solution to
${\cal E}_{u,g,f}$, where $h:{\Bbb R}\to{\Bbb R}$ is an arbitrary
smooth  function.}
\paragraph{3.4 Example.}
As an illustration of Theorem 3.3, consider the equation
\be u_t+u.u_x=e^u. \label{eq:p2}\ee
Then $g=u$, $f=e^u$, $\varphi =-e^{-u}$, and $\ell=u.(1-\ln(-u))$.
Therefore, by the Theorem 3.3, the solution of equation
(\ref{eq:p2}) is $x=(u+1).e^{-u}+h(t+e^{-u})$, where $h:{\Bbb
R}\to{\Bbb R}$ is an arbitrary function. For example, for
$h(s)=s$, we have the solution $u(x,t)=-\textrm{LW}(x-t)$.
\paragraph{3.5 Example.}
As an another illustration, consider the equation
\be u_t+(u^m)_x=u^n, \label{eq:p3}\ee
where $m$ and $n$ are integers. If $1<n\neq m$. Then
$g=m.u^{m-1}$, $f=u^n$, and
$$\ell(u)=\left\{ \begin{array}{lcl}
\di \frac{m}{m-n}.\big((1-n).u\big)^{(n-m)/(n-1)} && \mbox{if}\;
1<n\neq m,
\\[3mm]
\di \frac{m}{1-m}.\ln\big((1-m).u\big) && \mbox{if}\; n=m\neq1, \\[3mm]
\di \frac{m}{m-1}.\exp\big((m-1).u\big) && \mbox{if}\; n=1\;
\mbox{and}\;m\neq1,
\\[3mm]
\di e^u && \mbox{if}\;n=m=1,
\end{array} \right. $$
For example, $u(x,t)=x.(1+e^{-t})/2$ is a solution to
$u_t+(u^2)_x=u$; and
$$u(x,t)=\frac{1}{6}\Big(x-3.t\pm\sqrt{x^2+9.t^2-6.x.t-36}\;\Big)$$
is a solution to $u_t+(u^3)_x=u^2$.


\begin{thebibliography}{9}
%
\bibitem{[5]} {\sc
Burgers, J.M.}, {\em The nonlinear diffusion equation}, Reidel,
Dordrecht, 1974.
\bibitem{[1]} {\sc Ibragimov, N.H.} (Editor), {\em CRC Handbook of Lie Group Analysis
of Differential Equations}, Vol. 1, Symmetries, Exact Solutions
and Conservation Laws, CRC Press, Boca Raton, 1994.
\bibitem{[3]} {\sc Olver, P.J.}, {\em Application of Lie
groups to Differential Equation}, Springer, New York, 1986.
\bibitem{[2]} {\sc Polyanin, A.D.} and {\sc Zaitsev, V.F.}, {\em Handbook of Nonlinear
Partial Differential Equations}, Chapman \& Hall/CRC, Boca Raton,
2004.
\bibitem{[4]} {\sc Smaoui, N } and {\sc Mekkaoui, M.}, {\em The generalized Burgers equation with and without a time
delay}, Journal of Applied Mathematics and Stochastic Analysis
Volume 2004 (2004), 1, pp. 73-96.
\bibitem{[6]} {\sc Tychynin, V.} and {\sc Pasin, O.}, {\em Nonlocal Symmetry and Generating Solutions
for the Inhomogeneous Burgers Equation}, Proceedings of Institute
of Mathematics of NAS of Ukraine 2004, Vol. 50, Part 1, 277-281.
\bibitem{[7]} {\sc Wood, L.W.}, {\em An exact solution for Burger's
equation}, Commun. Numer. Meth. Engng 2006; 22: 797-798.
%
\end{thebibliography}
\end{document}